\documentclass{svproc}
\usepackage{amssymb}
\usepackage{amsmath}
\usepackage{graphicx}
\usepackage[cmtip,all]{xy}

\usepackage{url}

\newcommand{\cP}{\mathcal{P}}

\newcommand{\dP}{\Vec\cP}

\newcommand\blfootnote[1]{%
  \begingroup
  \renewcommand\thefootnote{}\footnote{#1}%
  \addtocounter{footnote}{-1}%
  \endgroup
}

\begin{document}
\mainmatter              
\title{On the Homomorphism Order of Oriented Paths and Trees.}
\titlerunning{Homomorphisms Order of Paths and Trees}  
%
\author{J. Hubi\v cka\inst{1} \and J. Ne\v set\v ril\inst{2} \and P. Oviedo \inst{3} \and O. Serra\inst{4}}
\authorrunning{Hubi\v cka et al.} 
%
%
\institute{Charles University, Department of Applied Mathematics (KAM), Praha 1, Czech Republic,\\
\email{hubicka@kam.mff.cuni.cz}
\and
Charles University, Computer Science Institute (IUUK), Praha 1, Czech Republic,\\
\email{nesetril@iuuk.mff.cuni.cz}
\and
University of Birmingham, School of Mathematics, Birmingham, UK,\\
\email{pxo006@student.bham.ac.uk}
\and
Universitat Politècnica de Catalunya, Department of Mathematics, Barcelona, Spain,\\
\email{oriol.serra@upc.edu}
}

\maketitle              

\begin{abstract}
A partial order is universal if it contains every countable partial order as a suborder. In 2017, Fiala, Hubi\v cka, Long and Ne\v set\v ril showed that every interval in the homomorphism order of graphs is universal, with the only exception being the trivial gap $[K_1,K_2]$. We consider the homomorphism order restricted to the class of oriented paths and trees. We show that every interval between two oriented paths or oriented trees of height at least 4 is universal. The exceptional intervals coincide for oriented paths and trees and are contained in the class of oriented paths of height at most 3, which forms a chain.
\end{abstract}

\blfootnote{The first and second author were supported by the project 21-10775S of 
the  Czech  Science Foundation (GA\v CR). This is part of a project that 
has received funding from the European Research Council (ERC) under the EU 
Horizon 2020 research and innovation programme (grant agreement No 
810115). The fourth author was supported by the Spanish Research Agency under project MTM2017-82166-P.
}

\section{Introduction}\label{section1}

Let $G_1,G_2$ be two finite directed graphs. A \emph{homomorphism} $f:G_1\to G_2$ is an arc preserving map $f:V(G_1)\to V(G_2)$. If such a map exists we write $G_1\leq G_2$. The relation $\leq$ defines a quasiorder on the class of directed graphs which, by considering homomorphic equivalence classes, becomes a partial order. A \emph{core} of a digraph is its minimal-size homomorphic equivalent digraph. As the core of a digraph is unique up to isomorphism, see \cite{llibre}, it is natural to choose a core as the representative of each homomorphic equivalence class.

Given a partial order $(\cP,\leq)$ and $a,b\in \cP$ satisfying $a<b$, the interval $[a,b]$ is a \emph{gap} if there is no $c\in \cP$ such that $a<c<b$. We say that a partial order is \emph{dense} if it contains no gaps. Finally, a partial order is \emph{universal} if it contains every countable partial order as a suborder.

The homomorphism order of graphs has proven to have a very rich structure \cite{llibre}. For instance, Welzl showed that undirected graphs, except for the gap $[K_1,K_2]$, are dense \cite{origdens}. Later, Ne\v set\v ril and Zhu \cite{pathhomomorphism} proved a density theorem for the class of oriented paths of height at least 4. Recently, Fiala et al. \cite{fractal} showed that every interval in the homomorphism order of undirected graphs, with the only exception of the gap $[K_1,K_2]$, is universal. The question whether this ``fractal" property was also present in the homomorphism order of other types of graphs was formulated. In this context, the following result was shown when considering the homomorphism order of the class oriented trees.

\begin{theorem}[\cite{tfg}]\label{theorem:proper_tree}
Let $T_1,T_2$ be oriented trees such that $T_1<T_2$. If the core of $T_2$ is not a path, then the interval $[T_1,T_2]$ is universal.
\end{theorem}

Theorem \ref{theorem:proper_tree} was presented in the previous Eurocomb \cite{tfg} but the result was not published as intervals of the form $[T_1,T_2]$ where the core of $T_2$ is a path remained to be characterized. In particular, the characterization of the universal intervals in the class of oriented paths was left as an open question. Here, we answer these questions and complete the characterization of the universal intervals in the class of oriented paths and trees by proving the following results (see Section 2 for definitions and notation).

\begin{theorem}\label{theorem:oriented_paths}
Let $P_1,P_2$ be oriented paths such that $P_1<P_2$. If the height of $P_2$ is greater or equal to 4 then the interval $[P_1,P_2]$ is universal. If the height of $P_2$ is less or equal to 3 then the interval $[P_1,P_2]$ forms a chain.
\end{theorem}

\begin{theorem}\label{theorem:oriented_trees}
Let $T_1,T_2$ be oriented trees such that $T_1<T_2$. If the height of $T_2$ is greater or equal to 4 then the interval $[T_1,T_2]$ is universal. If the height of $T_2$ is less or equal to 3 then the interval $[T_1,T_2]$ forms a chain.
\end{theorem}

It is known that the core of an oriented tree of height at most 3 is a path. Thus, oriented paths of height less or equal to 3 are the only exception when considering the presence of universal intervals in both the class of oriented paths and the class of oriented trees. Hence, the nature of its intervals in terms of density and universality is completely determined.

\begin{figure}[ht]
\centering
\includegraphics[scale=0.6]{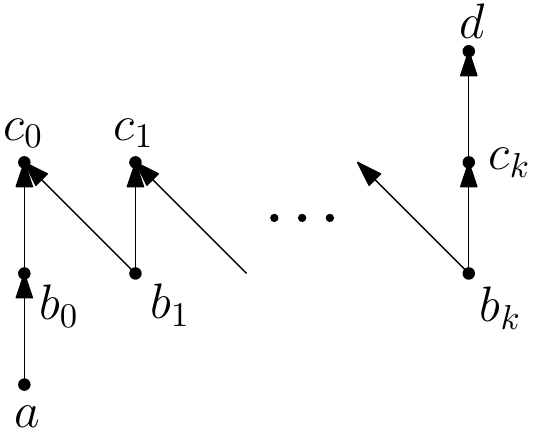}
\caption{The path $L_k$.}
\label{fig.L}
\end{figure}

Let $\Vec{P}_n$ denote the directed path consisting on $n+1$ vertices and $n$ consecutive forward arcs. Let $L_k$ denote the oriented path with vertex sequence $(a,b_0,c_0,b_1,c_1,\dots,b_k,c_k,d)$ and arcs $(ab_0,b_0c_0,b_1c_0,b_1c_1,\dots,b_kc_{k-1},$ $b_kc_k,c_kd)$, as shown in Figure \ref{fig.L}.

It is easy to check that $\Vec{P_0}<\Vec{P_1}<\Vec{P_2}$ and there is no oriented path strictly between $P_0$ and $P_1$ nor $P_1$ and $P_2$. Among oriented paths of height equal to three, the only cores are the paths $L_k$ for all $k\geq 0$, and it is also easy to see that $L_k\leq L_m$ if and only if $k\geq m$. Combining these observations we get that oriented paths of height at most 3 are found at the ``bottom" of the homomorphism order of paths and trees and form the following chain
$$\Vec{P_0}<\Vec{P_1}<\Vec{P_2}<\dots<L_{k+1}<L_k<L_{k-1}<\dots<L_2<L_1<L_0=\Vec{P_3}.$$

From the above results we can deduce the following corollary for both the homomorphism order of paths and the homomorphism order of trees.

\begin{corollary}\label{corollary}
An interval $[G_1,G_2]$ in the homomorphism order of oriented paths (resp. trees) is universal if and only if it contains two paths (resp. trees) which are incomparable in the homomorphism order.
\end{corollary}

It might seem that combining Theorems 1 and 2 one would get a proof of Theorem 3. However, this is not true as we have to consider intervals of form $[T_1,P_2]$ where the core of $T_1$ is not a path and the height of $P_2$ is greater or equal to 4. For this reason, we need a new density theorem.

\begin{theorem}\label{theorem:density}
Let $T_1$ be an oriented tree and $P_2$ a oriented path such that $T_1<P_2$. If the height of $P_2$ is greater or equal to 4 then there exists an oriented tree $T$ satisfying $T_1<T<P_2$.
\end{theorem}

The proof of Theorem \ref{theorem:density} is a more general version of the density theorem for paths proved by Ne\v set\v ril and Zhu \cite{pathhomomorphism}. Due to space limitations we omit the proof in this note.

It is now straightforward to prove Theorem 3. Given oriented trees $T_1,T_2$ with $T_1<T_2$ and $T_2$ with height greater or equal to 4, we consider two cases: if the core of $T_2$ is not a path we are done by applying Theorem \ref{theorem:proper_tree}; else, the core of $T_2$ is a path so by Theorem \ref{theorem:density} there exists a tree $T$ satisfying $T_1<T<T_2$. Now, if the core of $T$ is a path we can apply Theorem \ref{theorem:oriented_paths} to the interval $[T,T_2]$, else the core of $T$ is not a path and we can apply Theorem \ref{theorem:proper_tree} to the interval $[T_1,T]$.


\section{Notation}\label{section:notation}

A \emph{(oriented) path} $P=(V(P),A(P))$ is a sequence of \emph{vertices} $V(P)=(p_0,\ldots ,p_n)$ together with a sequence of \emph{arcs} $A(P)=(a_1,\dots,a_n)$ such that, for each $1\le i\leq n$, $a_i=p_{i-1}p_i$ or $p_{i-1}p_i$. We denote by $i(P)=p_0$ and $t(P)=p_n$ the \emph{initial} and \emph{terminal} vertex of $P$ respectively. We denote by $\dP$ the set of all oriented paths. Let $\Vec P_n$ denote the path with vertex set $V(\Vec{P_{n}})=(0,1,\dots,n)$ and arc set $A(\Vec{P_n})=(01,12,\dots,(n-1)n)$. The \emph{height} of a path $P$, denoted $h(P)$, is the minimum $k\geq 0$ such that $P\to \Vec{P_k}$. It is known that for any path $P$ there is a unique homomorphism $l:P\to \Vec{P}_{h(P)}$. Given such homomorphism $l:P\to \Vec{P}_{h(P)}$, the \emph{level} of a vertex $v\in V(P)$ is the integer $l(v)$. We define the \emph{level} of an arc $a\in A(P)$ as the greatest level of its incident vertices and denote it $l(a)$.

We write $P^{-1}$ for the \emph{reverse} path of $P$, where $V(P^{-1})=(p_n,p_{n-1},\ldots ,p_0)$ and $pp'\in A(P^{-1})$ if and only if $p'p\in A(P)$. That is, $P^{-1}$ is the path obtained from flipping the path $P$. 
Given paths $P,P'$, the \emph{concatenation} $PP'$ is the path with vertex set $V(PP')=(p_0,\ldots ,p_n=p'_0,p'_1,\ldots ,p'_{n'})$ and arc set $A(PP')=(a_1,\dots,a_n,a'_1,\dots,a'_{n'})$.

A \emph{zig-zag of length} $n$ consists on a path $Z$ with $n+1$ vertices whose arcs consecutively alternate from forward to backward. Note that all arcs in a zig-zag must have the same level. When considering a zig-zag $Z$ as a subpath of some path $P$, we define the \emph{level of Z} as the level of its arcs in $P$, and denote it $l(Z)$.


\section{Proofs}\label{section:proof}

The proof of Theorem \ref{theorem:oriented_paths} consists of the construction of an embedding from the homomorphism order of all oriented paths, which was proven to be universal \cite{univepaths}, into the interval $[P_1,P_2]$. The construction is based on the standard indicator technique initiated by Hedrlín and Pultr \cite{indicator}, see also \cite{llibre}. The method takes oriented paths $Q$ and $I$ and creates the path $\Phi_I(Q):= Q\ast I(a,b)$, obtained from $Q$ by replacing each arc $qq'\in A(Q)$ with a copy of the path $I$ by identifying the vertex $q$ with $i(I)$ and the vertex $q'$ with $t(I)$. If the path $I$ is well chosen we can guarantee that $\Phi_I$ induces an embedding from $(\dP,\leq)$ into the interval $[P_1,P_2]$. For an arc $qq'\in Q$, we denote $I_{qq'}:=\Phi_I(qq')$ to the copy of $I$ replacing $qq'$.

\begin{lemma}\label{lemma:embedding_paths}
Let $I$ be a path. Let $I_1,I_2$ be copies of $I$ and let $\epsilon_1,\epsilon_2\in \{-1,1\}$. Suppose that, for every $Q\in \dP$, the following conditions hold:
\begin{enumerate}
\item[(i)] $P_1<\Phi_I(Q)<P_2$;
\item[(ii)] Every homomorphism $f:I\to \Phi_I(Q)$ satisfies $f(I)\subseteq I_{qq'}$ for some $qq'\in A(Q)$;
\item[(iii)] Every homomorphism $g:I_1^{\epsilon_1}I_2^{\epsilon_2}\to \Phi_I(Q)$ satisfies that for every $z,z'\in \{`i\text{'} ,`t\text{'} \}$,  
$$z(I_1)=z'(I_2) \text{ implies } z(I_{q_1q_1'})=z'(I_{q_2q_2'})$$
where $g(I_1)\subseteq I_{q_1q_1'}$ and $g(I_2)\subseteq I_{q_2q_2'}$.
\end{enumerate}
Then $\Phi_I$ is a poset embedding from the class $\dP$ of paths into the interval $[P_1,P_2]$.
\end{lemma}

The power of the method is that the construction of the embedding reduces to finding a suitable gadget $I$ satisfying the above conditions.

To construct $I$, first consider a path $P$ such that $P_1<P<P_2$. We know such a path exist since paths of height greater or equal to 4 are dense \cite{pathhomomorphism}. Without loss of generality let $P$ be a core. Consider a surjective homomorphism $h:P\to P_2$. We can assume that such homomorphism always exists, since otherwise we will be able to find another path $P'$, with $P<P'<P_2$, which admits a surjective homomorphism into $P_2$. Then, since $h$ is surjective, there exist two different vertices $v_1,v_2\in V(P)$ such that $h(v_1)=h(v_2)$. Without loss of generality suppose that $v_1$ appears before $v_2$ in the sequence $V(P)$. Note that $l(v_1)=l(v_2)$ as homomorphisms preserve level vertex difference. The path $P$ is then naturally split into three subpaths: the subpath $A$ from $i(P)$ to $v_1$, the subpath $B$ from $v_1$ to $v_2$, and the subpath $C$ from $v_2$ to $t(P)$. So $P=ABC$. 

Let $a_1$ be the last arc of $A$ and let $a_2$ be the first arc of $C$. 

We first assume the case that neither $B=Z'_1$, $B=Z_2'$ nor $B=Z'_1Z'_2$ for some zig-zags $Z'_1$ and $Z'_2$ of level $l(a_1)$ and $l(a_2)$ in $P$ respectively.

Next, we choose as a gadget the following path $$I=Z_1A^{-1}ABCC^{-1}Z_2$$
where $Z_1$ and $Z_2$ are long enough zig-zags of even length and levels $l(a_1)$ and $l(a_2)$ in $I$ respectively.

It only remains to check that $I$ satisfies the conditions of Lemma \ref{lemma:embedding_paths}.

Let $Q\in \dP$. It follows by construction that $P_1<\Phi_I(Q)<P_2$. To see this, first note that $P\subset I$ so it is clear that $P_1<\Phi_I(Q)$. Now, observe that there is a homomorphism $\rho:I\to P$ which collapses the zig-zags $Z_1$ and $Z_2$ into the arcs $a_1$ and $a_2$ respectively, and maps $A^{-1}$ into $A$ and $C^{-1}$ into $C$ identifying the correspondent vertices. Then, the map $\rho':\Phi_I(Q)\to P_2$ defined for each copy $I_{qq'}$ of $I$ as $\rho'(I_{qq'}):=(h\circ\rho)(I_{qq'})$ is well defined since $h(v_1)=h(v_2)$. Finally, suppose that there is a homomorphism $f:P_2\to \Phi_I(Q)$. Since the zig-zags are long enough we must have $f(P_2)\subset I_{qq'}$ for some $qq'\in A(Q)$, implying that $P_2\to I\to P$, a contradiction. Thus, $I$ satisfies condition \textit{(i)} of Lemma \ref{lemma:embedding_paths}.

We say that a digraph is \emph{rigid} if its only automorphism is the identity map. To check conditions \textit{(ii)} and \textit{(iii)} of Lemma \ref{lemma:embedding_paths} we shall use several times the fact that the core of a path is rigid, see Lemma 2.1 in \cite{pathhomomorphism}. 

Another trivial property of homomorphisms is also that they can never increase the distance between two vertices when considering its images, see \cite{llibre}. That is, every homomorphism $f:G_1\to G_2$ satisfies $d(v_1,v_2)\geq d(f(v_1),f(v_2))$ for every $v_1,v_2\in V(G_1)$.

Let $f:I\to \Phi_I(Q)$ be a homomorphism. Since the zig-zags $Z_1$ and $Z_2$ are long enough, we have, for $ABC\subset I$, that $f(ABC)\subset I_{qq'}$ for some $qq'\in A(Q)$. Because $I_{qq'}\to P$ and $P$ is rigid, we must have that $f$ maps the path $ABC$ in $I$ to the copy of $ABC$ in $I_{qq'}$ via the identity map. If follows by a distance argument that $f(I)\subset I_{qq'}$, so condition \textit{(ii)} holds.

Finally, let $g:I_1I_2\to \Phi_I(Q)$ be a homomorphism. Note that $t(I_1)=i(I_2)$. Let $q_1q_1',q_2q_2'\in A(Q)$ such that $g(I_1)\subset I_{q_1q_1'}$ and $g(I_2)\subset I_{q_2q_2'}$. Let $w_1$ be the terminal vertex of the subpath $ABC$ in $I_1$ and let $w_2$ be the initial vertex of the subpath $ABC$ in $I_2$. Recall from the paragraph above that any homomorphism $f:I\to \Phi_I(Q)$ maps the path $ABC$ in $I$ to the path $ABC$ in $I_{qq'}$ fixing all vertices. Thus, $g(w_1)$ is the terminal vertex of the subpath $ABC$ in $I_{q_1q_1'}$ and $g(w_2)$ is the initial vertex of the subpath $ABC$ in $I_{q_2q_2'}$. We can not have $q_1=q_2'$ since in such case $d(g(w_1),g(w_2))>d(w_1,w_2)$. Both cases $q_1=q_2$ and $q_1'=q_2'$ imply that there is a homomorphism $g':I_1I_2\to I$. Suppose that this is the case. Observe that $w_1$ is joined to $w_2$ by the path $C^{-1}Z_1Z_2A^{-1}$ in $I_1I_2$. We would also have that $g'(w_2)$ is joined to $g'(w_1)$ by the path $ABC$ in $I$. Thus, we must have $B=Z_1'$, $Z_2'$ or $Z'_1Z'_2$ for some zig-zags of level equal to $Z_1$ and $Z_2$ respectively, which contradicts our first assumption. The only remaining case is $q_1'=q_2$.

The above argument is valid for the case $\epsilon_1=\epsilon_2=1$ in condition \textit{(iii)}. However, the other cases are analogous if not even simpler. Hence, applying Lemma \ref{lemma:embedding_paths} to the path $I$ completes the proof.

We have assumed that neither $B= Z_1'$, $B=Z_2'$ nor $B=Z'_1Z'_2$. If this is not the case what we do is take $l(a_1)\neq l(a_2)$ and consider the auxiliary path
$$P'=ABCC^{-1}Z_2CC^{-1}Z_1A^{-1}ABC$$
if $B\neq Z_1'$, or
$$P'=ABCC^{-1}Z_2A^{-1}AZ_1A^{-1}ABC$$
otherwise, where again $Z_1$ and $Z_2$ are long enough zig-zags of even length and level $l(a_1)$ and $l(a_2)$ in $P'$ respectively. Then we show that $P'$ satisfies the required properties for $P$ above and we apply the arguments to $P'$.


\section{Final Remarks}

For directed graphs, and even more, for general relational structures, the gaps of the homomorphism order are characterized by Ne\v set\v ril and Tardif in full generality \cite{dual}. They are all related to trees. However the homomorphism order of these general relational structures does not enjoy the simplicity of the ``bottom" of the homomorphism order of paths and trees. We think that Corollary \ref{corollary} may shed light on the characterization of universal intervals in the homomorphism order of these more general structures.


\end{document}